\documentclass[]{amsart}
\def\crn#1#2{{\vcenter{\vbox{\hbox{\kern#2pt \vrule width.#2pt height#1pt}
                                                    \hrule height.#2pt}}}}

\newtheorem{thm}{Theorem}

\begin{document}
\def\thefootnote{}
\footnotetext{MGE was supported as a Senior Research Fellow of the Australian
Research Council.
PWM was supported by `Fonds zur F\"orderung der
wissenschaftlichen Forschung, Projekt P~10037~PHY'.
The authors would also like to thank the organizers of the Nineteenth
Winter School on Geometry and Physics held in Srn\'\i, the Czech Republic, in
January 1999, where discussions concerning this article were initiated.\\
This paper is in final form and no version of it will be submitted for
publication elsewhere.}
\begin{center}{\LARGE\bf Some Remarks on the Pl\"ucker Relations}\\[15pt]
{\Large\bf Michael G. Eastwood and Peter W.\ Michor}
\end{center}

\section{The Pl\"ucker relations}

Let $V$ denote a finite-dimensional vector space. An $s$-vector
$P\in\Lambda^sV$ is called {\sl decomposable} or {\sl simple}
if it can be written
in the form
$$P=u\wedge v\wedge\cdots\wedge w\quad\mbox{for }u,v,\ldots,w\in V.$$
We shall use in the following both Penrose's abstract index notation
and exterior calculus with the conventions of \cite{G}.

\begin{thm}\label{plueckerthm}
Let $P\in\Lambda^sV$ be an $s$-vector. Then $P$ is decomposable if
and only if one of the following conditions holds:
\begin{enumerate}
\item\label{pluecker} $i(\Phi)P\wedge P=0$ for all
       $\Phi\in \Lambda^{s-1}V^*$.
       In index notation $P_{[abc\cdots d}P_{e]fg\cdots h}=0$.
\item\label{dualpluecker}
       $i(i_P\Psi)P=0$ for all $\Psi\in \Lambda^{s+1}V^*$.
\item\label{contraction} $i_{\alpha_1\wedge\dots\wedge\alpha_{s-k}}P$ is
       decomposable for all $\alpha_i\in V^*$, for any fixed
       $k\ge 2$.
\item\label{improvedpluecker}
       $i(\Psi)P\wedge P=0$ for all $\Psi\in \Lambda^{s-2}V^*$
       In index notation $P_{[abc\cdots d}P_{ef]g\cdots h}=0$.
\item\label{dualimprovedpluecker}
       $i(i_P\Psi)P=0$ for all $\Psi\in \Lambda^{s+2}V^*$.
\end{enumerate}
\end{thm}
\begin{proof}
(\ref{pluecker}) These are the well known classical Pl\"ucker relations.
For completeness' sake we include a proof.
Let $P\in\Lambda^n V$ and consider the induced linear mapping
$\sharp_P:\Lambda^{s-1}V^*\to V$. Its image, $W$, is contained in each linear
subspace $U$ of $V$ with $P\in\Lambda^s U$. Thus $W$ is the minimal subspace
with this property. $P$ is decomposable if and only if $\dim W= s$, and this
is the
case if and only if $w\wedge P=0$ for each $w\in W$. But $i_\Phi P$ for
$\phi\in \Lambda^{s-1}V^*$ is the typical element in $W$.

(\ref{dualpluecker}) This well known variant of the
Pl\"ucker relations follows by duality (see \cite{GH}):
\begin{align*}
\langle  P\wedge i(\Phi)P, \Psi \rangle &=
\langle  i(\Phi)P, i_P\Psi \rangle =
\langle  P, \Phi\wedge i_P\Psi \rangle =
\\
&=(-1)^{(s-1)}\langle  P, i_P\Psi\wedge \Phi \rangle =
(-1)^{(s-1)}\langle  i(i_P\Psi)P, \Phi \rangle.
\end{align*}

(\ref{contraction}) This is due to \cite{MV}. There it is proved
using exterior algebra.
Apparently, this result is included in formula (4), page 116
of \cite{W}.

(\ref{improvedpluecker}) Another proof using representation theory will
be given below.
Here we prove it by induction on $s$.
Let $s=3$. Suppose that $i_\alpha P\wedge P = 0$ for all $\alpha\in V^*$.
Then for all $\beta\in V^*$ we have
$ 0 = i_\beta(i_\alpha P\wedge P)
= i_{\alpha\wedge\beta}P\wedge P + i_\alpha P\wedge i_\beta P$.
Interchange $\alpha$ and $\beta$ in the last expression and add it to the
original, then we get
$0 = 2i_\alpha P\wedge i_\beta P$ and in turn
$i_{\alpha\wedge\beta}P\wedge P= 0$ for all $\alpha$ and $\beta$, which are
the original Pl\"ucker relations, so $P$ is decomposable.
Now the induction step.
Suppose that $P\in\Lambda^s V$ and that
$i_{\alpha_1\wedge\dots\wedge\alpha_{s-2}}P\wedge P = 0$
for all $\alpha_i\in V^*$.
Then we have
$$0= i_{\alpha_1}(i_{\alpha_1\wedge\dots\wedge\alpha_{s-2}}P\wedge P) =
i_{\alpha_1\wedge\dots\wedge\alpha_{s-2}}P\wedge i_{\alpha_1}P =
i_{\alpha_2\wedge\dots\wedge\alpha_{s-2}}
(i_{\alpha_1}P)\wedge (i_{\alpha_1}P)$$
for all $\alpha_i$, so that by induction we may conclude that
$i_{\alpha_1}P$ is decomposable for all $\alpha_1$, and then by
(\ref{contraction}) $P$ is decomposable.

(\ref{dualimprovedpluecker})
Again this follows by duality.
\end{proof}

Let us note that the following result
(Lemma 1 in \cite{GM}), a version of the `three plane lemma'
also implies (\ref{contraction}):

Let $\{P_i: i\in I\}$ be a family of decomposable non-zero
$k$-vectors in $V$ such that each $P_i+P_j$ is again
decomposable. Then
\begin{enumerate}
\item[(a)] either the linear span $W$ of the linear subspaces
       $W(P_i)=\operatorname{Im}(\sharp_{P_i})$ is at most
       $(k+1)$-dimensional
\item[(b)] or the intersection $\bigcap_{i\in I}W(P_i)$ is at least
       $(k-1)$-dimensional.
\end{enumerate}

Finally note that (\ref{pluecker}) and
(\ref{improvedpluecker}) are both invariant under ${\mathrm{GL}}(V)$.
In the next section we shall decompose (\ref{pluecker}) into its
irreducible components in this representation.

If $\dim V$ is high enough in comparison with $s$, then
(\ref{improvedpluecker}) seemingly comprises less equations.

\section{Representation theory}
In order efficiently to analyse (\ref{pluecker}) and
(\ref{improvedpluecker})
it is necessary to take a small excursion through representation
theory. An extensive discussion of Young tableau may be found
in~\cite{fulton}. Here we shall just need
$$Y^{s,t}\equiv\raisebox{-30pt}{\begin{picture}(60,80)
          \put(20,0){\line(1,0){10}}
          \put(20,10){\line(1,0){10}}
          \put(20,30){\line(1,0){10}}
          \put(20,40){\line(1,0){20}}
          \put(20,50){\line(1,0){20}}
          \put(20,70){\line(1,0){20}}
          \put(20,80){\line(1,0){20}}
          \put(20,0){\line(0,1){80}}
          \put(30,0){\line(0,1){80}}
          \put(40,40){\line(0,1){40}}
          \put(25,23){\makebox(0,0){$\vdots$}}
          \put(25,63){\makebox(0,0){$\vdots$}}
          \put(35,63){\makebox(0,0){$\vdots$}}
          \put(5,40){\makebox(0,0){$s$}}
          \put(55,60){\makebox(0,0){$t$}}
          \put(15,40){\makebox(0,0){$\left\{\rule{0pt}{44pt}\right.$}}
          \put(45,60){\makebox(0,0){$\left.\rule{0pt}{22pt}\right\}$}}
\end{picture}}$$
regarded as irreducible representations of ${\mathrm{GL}}(V)$. Then, as
special cases of the Littlewood-Richardson rules, we have
$$\begin{array}{ccr}
\Lambda^sV\otimes\Lambda^sV&=&Y^{s,s}\oplus Y^{s+1,s-1}\oplus
                              Y^{s+2,s-2}\oplus Y^{s+3,s-3}
                              \oplus\cdots\oplus Y^{2s,0}\\[3pt]
\Lambda^{s+1}\otimes\Lambda^{s-1}V&=&       Y^{s+1,s-1}\oplus
                              Y^{s+2,s-2}\oplus Y^{s+3,s-3}
                              \oplus\cdots\oplus Y^{2s,0}\\[3pt]
\Lambda^{s+2}\otimes\Lambda^{s-2}V&=&
                              Y^{s+2,s-2}\oplus Y^{s+3,s-3}
                              \oplus\cdots\oplus Y^{2s,0}
\end{array}$$
and from the first two of these (\ref{pluecker}) says that
$P\otimes P\in Y^{s,s}$. In fact,
$$(\star\star)\qquad\begin{array}[t]{cccccccc}
\Lambda^sV\odot\Lambda^sV &=&Y^{s,s}&\oplus&Y^{s+2,s-2}&\oplus&\cdots\\[3pt]
\Lambda^sV\wedge\Lambda^sV&=&&Y^{s+1,s-1}&\oplus&Y^{s+3,s-3}&\oplus&\cdots
\end{array}$$
so we can also see the equivalence of (\ref{pluecker}) and
(\ref{improvedpluecker}) without any calculation.
Having decomposed $\Lambda^sV\odot\Lambda^sV$ into irreducibles, it behoves one
to investigate the consequences of having each irreducible component of
$P\otimes P$ vanish separately. The first of these gives us another improvement
on the classical Pl\"ucker relations:
\begin{thm}\label{optimal}
An $s$-form $P$ is simple if and only if the component of
$P\otimes P$ in $Y^{s+2,s-2}$ vanishes.\end{thm}
\begin{proof}The representation $Y^{s+2,s-2}$ may be realised as those tensors
$$T_{a_1b_1a_2b_2\ldots a_{s-2}b_{s-2}cdef}$$
which are symmetric in the pairs $a_jb_j$ for $j=1,2,\ldots s-2$, skew in
$cdef$, and have the property that symmetrising over any three indices gives
zero. The corresponding Young projection of
$$P_{a_1a_2\ldots a_{s-2}cd}P_{b_1b_2\ldots b_{s-2}ef}$$
is obtained by skewing over $cdef$ and symmetrising over each of the pairs
$a_jb_j$ for $j=1,2,\ldots,s-2$. Its vanishing, therefore, is equivalent to
the vanishing of
$$Q_{[cd}Q_{ef]}\quad\mbox{where }Q_{cd}=
  \alpha^{a_1}\beta^{a_2}\cdots\gamma^{a_{s-2}}P_{a_1a_2\ldots a_{s-2}cd}$$
for all $\alpha^a,\beta^a,\ldots,\gamma^a\in V^\ast$. According
to~(\ref{improvedpluecker}), this means that $Q_{cd}$ is simple. Therefore,
the theorem is equivalent to criterion (\ref{contraction}) of
Theorem~\ref{plueckerthm}.
\end{proof}

Notice that this generally cuts down further the number of equations needed
to characterise the simple $s$-vectors. The simplest instance of this is
for 4-forms:
$P$~is simple if and only if
$$P_{[abcd}P_{ef]gh}=P_{[abcd}P_{efgh]}.$$
Written in this way, it is slightly surprising that one can deduce the
vanishing of each side of this equation separately. Theorem~\ref{optimal}
is optimal in the sense that the vanishing of any other component or
components in the irreducible decomposition ($\star\star$) of $P\otimes P$
is either insufficient to force simplicity or causes $P$ to vanish. In the
case of four-forms, for example,
$$P_{[abcd}P_{efgh]}=0$$
if
$P=v\wedge Q$
for some vector $v$ and three-form~$Q$. On the other hand, if the $Y^{4,4}$
component of $P\otimes P$ vanishes, then arguing as in the proof of
Theorem~\ref{optimal} shows that $P=0$.

{\small P. W. Michor:
Institut f\"ur Mathematik, Universit\"at Wien,
Strudlhofgasse 4, A-1090 Wien, Austria, {\it and}:
Erwin Schr\"odinger Institute, Boltzmanngasse 9, A-1090, Wien,
Austria.} \\{\small E-mail: michor@esi.ac.at}\\
{\small M.\ Eastwood: Dept.\ Pure Math., Univ. of Adelaide,
Adelaide, SA 5005, Australia.} \\{\small E-mail:
meastwoo@maths.adelaide.edu.au}\\
\end{document}